\documentclass[11pt]{article}
\usepackage{latexsym}
\usepackage{amssymb}
\usepackage{amsmath, amsfonts, amscd, amsthm, mathrsfs}
\usepackage{mathptm}    
\usepackage[all,dvips,arc,curve,color,frame]{xy}

\title{A short guide to $p$-torsion of abelian varieties in characteristic $p$}
\author{Rachel Pries
\footnote{The author was partially supported by NSF grant DMS-04-00461.}}
\date{}

\newtheorem{theorem}{Theorem}[section]
\newtheorem{lemma}[theorem]{Lemma}
\newtheorem{proposition}[theorem]{Proposition}

\newenvironment{prf}[1]{\trivlist
\item[\hskip
\labelsep{\it #1.\hspace*{.3em}}]}{
\endtrivlist}

\newtheorem{predefinition}[theorem]{Definition}

\newtheorem{preremark}[theorem]{Remark}

\newtheorem{prenotation}[theorem]{Notation}

\newtheorem{preexample}[theorem]{Example}
\newenvironment{example}{\begin{preexample}\rm}{\end{preexample}}
\newtheorem{preclaim}[theorem]{Claim}

\newtheorem{prequestion}[theorem]{Question}

\makeatletter 
\def\emppsubsection{\@startsection{subsection}{2}{\z@}{-3.25ex plus -1ex minus -.2ex}{-1em}{\bf}}

\makeatother

\newcommand \CA {{\cal A}}

\newcommand \ZZ {{\mathbb Z}}
\newcommand  \FF {{\mathbb F}}

\newcommand \GG {{\mathbb G}}
\newcommand  \QQ {{\mathbb Q}}
\newcommand  \NN {{\mathbb N}}

\newcommand \Spec {\mathop{\rm Spec}}

\newcommand \Hom {\mathop{\rm Hom}}

\newcommand \dime {\mathop{\rm dim}}

\newcommand \codim {\mathop{\rm codim}}

\begin{document}
\maketitle

\begin{abstract}
There are many equivalent ways to describe the $p$-torsion of a
principally polarized abelian variety in characteristic $p$.
We briefly explain these methods and then illustrate them for abelian varieties $A$ of arbitrary dimension $g$
in several important cases, including when $A$ has $p$-rank $f$ and $a$-number $1$ 
and when $A$ has $p$-rank $f$ and $a$-number $g-f$.
We provide complete tables for abelian varieties of dimension up to four.
\end{abstract}

\section{Introduction}

In recent years, there have been many important results about the $p$-torsion of a
principally polarized abelian variety in characteristic $p$.
This $p$-torsion can be described in terms of a group scheme or a Dieudonn\'e module.  
It can be classified using its final type or its Young type. 
It can be identified with an element in the Weyl group of the sympletic group 
or with a cycle class in the tautological ring of $\CA_g$.

In this paper, we briefly summarize the main types of classification.  
We give a thorough description of the $p$-torsion of a principally polarized abelian variety $A$ of arbitrary dimension $g$
in several important cases, including when $A$ has $p$-rank $f$ and $a$-number $1$, and when $A$ has $p$-rank $f$ 
and $a$-number $g-f$.
We provide complete tables for the $p$-torsion types that occur for $g \leq 4$, including 
the sixteen types of $p$-torsion that occur for abelian varieties of dimension four.   
We hope that this paper will inspire the reader to learn more about the outstanding research in this area.

\section{Methods to classify the $p$-torsion}\label{S2}

Let $k$ be an algebraically closed field of characteristic $p$.
Let $\CA_g:=\CA_g \otimes \FF_p$ be the moduli space of
principally polarized abelian varieties of dimension $g$ defined over $k$.
For an abelian variety $A \in \CA_g(k)$, let $A[p]$ denote its $p$-torsion.
We summarize several different ways of describing $A[p]$. 

\subsection{Group schemes}

Let $A$ be an abelian variety of dimension $g$ defined over $k$.
The $p$-torsion $A[p]$ is a finite commutative group scheme annihilated by $p$ with rank $p^{2g}$
having homomorphisms $F$ (Frobenius) and $V$ (Vershiebung).
If $A$ is principally polarized, then ${\rm im}(F)={\rm ker}(V)$ and ${\rm im}(V)={\rm ker}(F)$.
Then $A[p]$ is called a quasi-polarized $BT_1$ $k$-group scheme (short for truncated Barsotti-Tate group of level 1).  The quasi-polarization implies that $A[p]$ is symmetric.
These group schemes were classified independently by Kraft (unpublished) \cite{Kraft} and by Oort \cite{O:strat}.
A complete description of this topic can be found in \cite{O:strat} or \cite{M:group}.

\begin{example} \label{Eord}
Let $\ZZ/p$ be the constant group scheme and let $\mu_p$ be the kernel of Frobenius on $\GG_m$. 
As a $k$-scheme, $\mu_p \simeq \Spec(k[x]/(x^p-1))$.
If $E$ is an ordinary elliptic curve then $E[p] \simeq \ZZ/p \oplus \mu_p$.
We denote this group scheme by $L$.
\end{example}

\begin{example} \label{EI1}
Let $\alpha_p$ be the kernel of Frobenius on $\GG_a$. 
As a $k$-scheme, $\alpha_p \simeq \Spec(k[x]/x^p)$.
The isomorphism type of the $p$-torsion of any two supersingular elliptic curves is the same.
If $E$ is a supersingular elliptic curve, we denote the isomorphism type of its $p$-torsion by $I_{1,1}$.
By \cite[Ex.\ A.3.14]{G:book}, $I_{1,1}$ fits into a non-split exact sequence of the form 
$0 \to \alpha_p \to I_{1,1} \to \alpha_p \to 0$.
The image of the embedded $\alpha_p$ is unique and is the kernel of both Frobenius and Verschiebung.
\end{example}
 
\begin{example} \label{EI2}
Let $A$ be a supersingular non-superspecial abelian surface.  
In other words, $A$ is isogenous, but not isomorphic, to the direct sum of two 
supersingular elliptic curves.
Let $I_{2,1}$ denote the isomorphism class of the group scheme $A[p]$. 
By \cite[Ex.\ A.3.15]{G:book}, there is a filtration
$H_1 \subset H_2 \subset I_{2,1}$ where $H_1 \simeq \alpha_p$,
$H_2/H_1 \simeq \alpha_p \oplus \alpha_p$, and $I_{2,1}/H_2 \simeq \alpha_p$.
If $G_1$ (resp.\ $G_2$) is the kernel of Frobenius (resp.\ Verschiebung)
then $G_1 \subset H_2$ and $G_2 \subset H_2$.
There is an exact sequence $0 \to H_1 \to G_1 \oplus G_2 \to H_2 \to 0$.
\end{example}

Two invariants of (the $p$-torsion of) an abelian variety are the $p$-rank and $a$-number.
The {\it $p$-rank} of $A$ is $f=\dime_{\FF_p} \Hom(\mu_p, A[p])$.
Then $p^f$ is the cardinality of $A[p](k)$.
The {\it $a$-number} of $A$ is $a=\dime_k \Hom(\alpha_p, A[p])$.
It is well-known that $0 \leq f \leq g$ and $0 \leq a \leq g-f$.
In Example \ref{Eord}, $f=1$ and $a=0$.
In Example \ref{EI1},  $f=0$ and $a=1$.
The group scheme $I_{2,1}$ in Example \ref{EI2} has $p$-rank $0$ 
since it is an iterated extension of copies of $\alpha_p$
and has $a$-number $1$ since ${\rm ker}(V^2)=G_1 \oplus G_2$ has rank $p^3$.

\subsection{Covariant Dieudonn\'e modules}

One can describe the $p$-torsion $A[p]$ using the theory of covariant Dieudonn\'e modules.
This is the dual of the contravariant theory found in \cite{Demazure}; see also \cite[A.5]{G:book}.
Briefly, let $\sigma$ denote the Frobenius automorphism of $k$. 
Consider the non-commutative ring ${\mathbb E}=k[F,V]$ with the relations 
$FV=VF=0$ and $F\lambda=\lambda^\sigma F$ and $\lambda V=V \lambda^\sigma$ for all $\lambda \in k$.
Let $(A,B)_\ell$ denote the left ideal ${\mathbb E}A+{\mathbb E}B$ of ${\mathbb E}$ generated by $A$ and $B$.
The Dieudonn\'e functor $D$ gives an equivalence of categories between $BT_1$ group schemes 
${\mathbb G}$ (with rank $p^{2g}$) and finite left ${\mathbb E}$-modules $D({\mathbb G})$
(having dimension $2g$ as a $k$-vector space).
If ${\mathbb G}$ is quasi-polarized, then there is a sympletic form on $D({\mathbb G})$.

For example, $D(\ZZ/p \oplus \mu_p) \simeq {\mathbb E}/(F, 1-V)_\ell \oplus {\mathbb E}/(V, 1-F)_\ell$, \cite[Ex.\ A.5.1 \& 5.3]{G:book}.
The Dieudonn\'e module for $I_{1,1}$ and $I_{2,1}$ can be found in Lemma \ref{Lunique}.

The $p$-rank of $A[p]$ is the dimension of $V^g D({\mathbb G})$.
The $a$-number of $A[p]$ equals $g-{\rm dim}(V^2D({\mathbb G}))$.

\subsection{Final types}
The isomorphism type of a symmetric $BT_1$ group scheme ${\mathbb G}$ over $k$
can be encapsulated into combinatorial data.
This topic can be found in \cite{O:strat}.
If ${\mathbb G}$ has rank $p^{2g}$, 
then there is a {\it final filtration} $N_1 \subset N_2 \subset \cdots \subset N_{2g}$ 
of $D({\mathbb G})$ as a $k$-vector space
which is stable under the action of $V$ and $F^{-1}$ so that $i={\rm dim}(N_i)$.
If ${\mathbb G}$ is quasi-polarized, then $N_{2g-i}$ and $N_i$ are orthogonal under the sympletic pairing.

The {\it final type} of ${\mathbb G}$ is $\nu=[\nu_1, \ldots, \nu_r]$ where ${\nu_i}={\rm dim}(V(N_i))$.
The final type of ${\mathbb G}$ is canonical, even if the final filtration is not.
There is a restriction $\nu_i \leq \nu_{i+1} \leq \nu_i +1$ on the final type.
All sequences satisfying this restriction occur.   
This implies that there are $2^g$ isomorphism types of symmetric $BT_1$ group schemes of rank $p^{2g}$.    
The $p$-rank is ${\rm max}\{i \ | \ \nu_i=i\}$ and the $a$-number is $g-\nu_g$.

Together with Ekedahl, Oort used this classification by final type to stratify $\CA_g$.  
The stratum of $\CA_g$ whose points have final type $\nu$ is locally closed and quasi-affine with 
dimension $\sum_{i=1}^g \nu_i$, \cite[Thm.\ 1.2]{O:strat}. 

\subsection{Young types}

Another combinatorial method to describe the isomorphism type of ${\mathbb G}$ uses a Young diagram.
This method was introduced by Van der Geer \cite{V:cycles} as a means of describing 
the Ekedahl-Oort strata in terms of degeneration loci for maps between flag varieties.

Given a final type $\nu$, let $\mu_j=\#\{i \ | \ 1 \leq i \leq g \ | \ j \leq i-\nu_i \}$.
Consider the Young diagram with $\mu_j$ squares in the $j$th row.
The {\it Young type} of ${\mathbb G}$ is $\mu=\{\mu_1, \mu_2, \ldots\}$.
The $p$-rank is $g-\mu_1$ and the $a$-number is $a={\rm max}\{j \ | \ \mu_j \not = 0\}$.
The codimension in $\CA_g$ of the stratum whose points have Young type $\mu$ is $\sum_{j=1}^a \mu_j$. 

\subsection{Elements of the Weyl group}

One can associate to $\mu$ 
an element $\omega$ of the Weyl group $W_g$ of the sympletic group $Sp_{2g}$, \cite{V:cycles}.
Here $W_g$ is identified with the subgroup of all $\omega \in S_{2g}$ so that $\omega(i)+\omega(2g+1-i)=2g+1$ for $1 \leq i \leq g$.
This subgroup is generated by the following involutions:
$s_i=(i,i+1)(2g-i,2g+1-i)$ for $1 \leq i <g$; and $s_g=(g,g+1)$. 

Given a Young type $\mu$, one defines $\omega$ as follows.
For $1 \leq i \leq g$, let $\omega(i)=c$ (respectively $\omega(i)=g+c$) 
if $i$ is the $c$th number such that $\mu_i=\mu_{i+1}$ 
(respectively $\mu_i \not =\mu_{i+1}$).  
For $1 \leq i \leq g$, let $\omega(2g+1-i)=2g+1-\omega(i)$.
This yields an element of $W_g$.  
One can express $\omega$ as a word in the involutions $s_1, \ldots, s_g$ of $S_{2g}$, 
although this expression is not unique.    

For example, in the ordinary case where $\mu = \emptyset$, then $\omega$ is given by 
$\langle 1, \ldots, 2g\rangle \stackrel{\omega} \to \langle g+1, \ldots, 2g, 1, \ldots, g \rangle$.
In the superspecial case where $\mu=\{g, \ldots, 1\}$, then $\omega={\rm id}$.
Further examples with $g \leq 4$ are in Section \ref{Stables}.

We briefly explain the importance of the Weyl group characterization.  
There is a second filtration of $D(\GG)$ which is stable under the action of $F$ and $V^{-1}$, 
which we denote by $N'_1 \subset N'_2 \subset \cdots \subset N'_{2g}$.
Then $\omega$ measures the interaction between these two filtrations.

For example, when $\GG$ is ordinary ($f=g$) then $N_g \cap N'_g=0$.  
Informally speaking, this means that the intersection of ${\rm Im}(V)$ and 
(a twist under $\sigma$ of) ${\rm Im}(F)$ is trivial. 
When $\GG$ is superspecial ($a=g$), then ${\rm dim}(N_i \cap N'_g)=i$ for $1 \leq i \leq g$.  
Informally speaking, this implies that $N_i$ is contained in (a twist under $\sigma$ of) ${\rm Im}(F)$. 
In general, ${\rm dim}(N_i \cap N_g') \geq i-\nu_i$. 
The $a$-number is ${\dime}(VD(\GG)\cap FD(\GG))=\dime(N_g \cap N_g')$.

One can identify the closures of the Ekedahl-Oort strata with cycle classes in the 
tautological ring of $\CA_g$.  
Let $\lambda_i$ for $1 \leq i \leq g$ be the Chern classes of the Hodge bundle of $\CA_g$.
These classes generate the tautological subring of $CH^*_{\QQ}(\CA_g)$
and satisfy $(1+\lambda_1 + \cdots + \lambda_g)(1-\lambda_1 + \cdots +(-1)^g\lambda_g)=1$, \cite[Thm.\ 1.1]{V:cycles}.  

\section{Important examples}

\subsection{Abelian varieties with $p$-rank $f$}

Given $g$ and $f$ such that $0 \leq f \leq g$, let $V_{g,f}$ denote the stratum of $\CA_g$ 
whose points correspond to principally polarized abelian varieties $A$ of dimension $g$ with $f_A \leq f$. 
Every component of $V_{g,f}$ has codimension $g-f$, \cite{NO}.
In this section, we describe the $p$-torsion that occurs for the generic point(s) of $V_{g,f}$.
The generic point of $V_{g,g}=\CA_g$ has $p$-rank $g$, $a$-number $0$, and $A[p] \cong L^g$. 
Using a dimension count, one can show that the generic point of every component of $V_{g,f}$ has $a$-number $1$ when $f <g$.

\begin{lemma} \label{Lunique}
Let $r \in \NN$.
There is a unique symmetric $BT_1$ group scheme of rank $p^{2r}$ with $p$-rank 0 and $a$-number 1,
which we denote $I_{r,1}$.  
The covariant Dieudonn\'e module of $I_{r,1}$ is ${\mathbb E}/(F^{r}+V^{r})_\ell$.
\end{lemma}

\begin{proof}
Let $I_{r,1}$ be a symmetric $BT_1$ group scheme of rank $p^{2r}$ with $p$-rank $0$ and $a$-number $1$.
It is sufficient to show that the final type of $I_{r,1}$ is uniquely determined.
The $p$-rank $0$ condition implies that $V$ acts nilpotently on $D(I_{r,1})$, so $\nu_1=0$.
The $a$-number $1$ condition implies that $r-1$ is the dimension of $V^2D(I_{r,1})$, so $\nu_{r}=r-1$.
The restrictions on $\nu_i$ imply 
that there is a unique final type possible for $I_{r,1}$, namely $\nu=[0,1, \ldots, r-1]$.

Consider $D={\mathbb E}/(F^r+V^r)_\ell$.  Note that $F^{r+1}=0$ and $V^{r+1}=0$ on $D$.
Then $D$ is an ${\mathbb E}$-module with dimension $2r$ as a $k$-vector space.
It has basis $\{F, \ldots, F^{r}, 1, V, \ldots, V^{r-1}\}$.
Then $VD$ has basis $\{V, \ldots, V^{r-1}, F^r\}$ and
$V^2D$ has basis $\{V^2, \ldots, V^{r-1}, F^r\}$.  Thus $D$ has $a$-number $1$.   
Continuing, one sees that $V$ is nilpotent on $D$ and thus the $p$-rank of $D$ is $0$.
Thus $D$ must be the covariant Dieudonn\'e module corresponding to $I_{r,1}$.
\end{proof}


\begin{proposition} \label{PVgf}
Let $A \in \CA_g(k)$ be a principally polarized abelian variety of dimension $g$ with $p$-rank $f$ and $a$-number $1$.
Then $A[p] \simeq L^{f} \oplus I_{g-f,1}$.  
The covariant Dieudonn\'e module of $A[p]$ is 
$$D \simeq ({\mathbb E}/(F, 1-V)_\ell \oplus {\mathbb E}/(V, 1-F)_\ell)^f \oplus {\mathbb E}/(F^{g-f}-V^{g-f})_\ell.$$
The final type of $A[p]$ is $\nu=[1, \ldots, f, f, \ldots, g-1]$.
The Young type is $\mu=\{g-f\}$.  
\end{proposition}

\begin{proof}
The decomposition of $A[p]$ must include $f$ copies of $L$ along with a group scheme of
rank $p^{2(g-f)}$ with $p$-rank $0$ and $a$-number $1$.  By Lemma \ref{Lunique}, 
the only possibility for the latter is $I_{g-f,1}$. 
The statement about the Dieudonn\'e module follows immediately.
For the final type, note that $\nu_g=g-1$ since $A[p]$ has $a$-number $1$ 
and $\nu_{f}=f$ since $A[p]$ has $p$-rank $f$.  
The numerical restrictions on $\nu_i$ then imply that $\nu=[1, \ldots, f, f, \ldots, g-1]$.
The Young type follows by direct calculation.
\end{proof}

If $f <g$, one can show that the group scheme $L^{f} \oplus I_{g-f,1}$ corresponds to the element $\omega$
of the Weyl group so that $\omega(f+1)=1$ and   
$\omega:\{1, \ldots, g\}-\{f+1\} \mapsto \{g+1, \ldots, 2g-1\}$ is increasing. 
The cycle class of the (reduced) stratum $V_{g,f}$ in the tautological ring of $\CA_g$ is given by 
$(p-1)(p^2-1)\ldots (p^{g-f}-1)\lambda_{g-f}$, \cite[Thm. 2.4]{V:cycles}.

\subsection{Abelian varieties with $a$-number $a$}

Given $g$ and $f$ such that $0 \leq f \leq g$, let $T_{g,a}$ denote the stratum of $\CA_g$ 
whose points correspond to principally polarized abelian varieties of dimension $g$ with $a_A \geq a$.
Then $T_{g,a}$ is irreducible unless $a=g$, \cite[Thm.\ 2.11]{V:cycles}.
In this section, we describe the $p$-torsion that occurs for the generic point(s) of $T_{g,a}$. 
It is well-known that $T_{g,a}$ has codimension $a(a+1)/2$.
The generic point(s) of $T_{g,a}$ have $a$-number $a$ and $p$-rank $g-a$.

\begin{proposition}
Let $A \in \CA_g(k)$ be an abelian variety of dimension $g$ with $p$-rank $f$ and $a$-number $g-f$.
Then $A[p] \simeq L^{f} \oplus (I_{1,1})^{g-f}$.
The covariant Dieudonn\'e module of $A[p]$ is 
$$D \simeq ({\mathbb E}/(F, 1-V)_\ell \oplus {\mathbb E}/(V, 1-F)_\ell)^f \oplus ({\mathbb E}/(F+V)_\ell)^{g-f}.$$
The final type is $\nu=[1, \ldots, f, \ldots, f]$.
The Young type is $\mu=\{g-f, g-f-1, \ldots, 1\}$ or $\emptyset$ if $g=f$.  
\end{proposition}

\begin{proof}
The decomposition of $A[p]$ must include $f$ copies of $L$ along with a group scheme of
rank $p^{2(g-f)}$ with $p$-rank $0$ and $a$-number $g-f$.  
The only possibility for the latter is $g-f$ copies of $D(I_{1,1})$. 
The statement about the Dieudonn\'e module follows immediately.
For the final type, note that $\nu_g=f$ since $A[p]$ has $a$-number $g-f$ 
and $\nu_{f}=f$ since $A[p]$ has $p$-rank $f$.  
The numerical restrictions on $\nu_i$ then imply that $\nu=[1, \ldots, f, \ldots, f]$.
The Young type follows by direct calculation.
\end{proof}

If $f >0$, one can show that the group scheme $L^{f} \oplus (I_{1,1})^{g-f}$ 
corresponds to the element $\omega$ of the Weyl group
$\langle 1, \ldots, 2g\rangle \stackrel{\omega} \to \langle g+1, \ldots, g+f, 1, \ldots a, g+f+1, \ldots 2g, a+1, \ldots, g\rangle$.
In \cite[Thm.\ 2.6]{V:cycles}, one finds a result on the cycle class of the (reduced) stratum $T_{g,a}$ 
in the tautological ring of $\CA_g$.

\subsection{Some indecomposable group schemes with $p$-rank $0$ and $a$-number $2$}

A symmetric $BT_1$ group scheme ${\mathbb G}$ is {\it indecomposable}
if ${\mathbb G} \not = {\mathbb G}_1 \oplus {\mathbb G}_2$ where 
${\mathbb G}_1$ and ${\mathbb G}_2$ are nontrivial symmetric $BT_1$ group schemes.
The group schemes $I_{r,1}$ are indecomposable.
We now describe an indecomposable group scheme $I_{r,2}$ of rank $p^{2r}$, $p$-rank $0$, and $a$-number $2$.

\begin{lemma} \label{a=2}  
Let $r \in \NN$ with $r \geq 3$.
Let $D={\mathbb E}/(F^{r-1}-V)_\ell \oplus {\mathbb E}/(V^{r-1}-F)_\ell$.
Then $D$ is the covariant Dieudonn\'e module of an indecomposable symmetric $BT_1$ group scheme 
with rank $p^{2r}$, $p$-rank $0$ and $a$-number $2$, which we   
denote by $I_{r,2}$.
It has final type $[0,1, \ldots, r-3,r-2,r-2]$ and Young type $\{r,1\}$.
\end{lemma}

\begin{proof}
The given decomposition of $D$ is the only possible decomposition of $D$ into covariant Dieudonn\'e modules, 
but neither of the factors in this decomposition is symmetric. 
Thus $I_{r,2}$ is indecomposable.  

Note that $F^r=0$ (resp.\ $V^r=0$) on the first (resp.\ second) factor of $D$.
Then $D=N_{2r}=\langle 1,F, \ldots, F^{r-1} \rangle \oplus \langle 1, V, \ldots, V^{r-1}\rangle$.
Thus $D$ has dimension $2r$ as a $k$-vector space and $I_{r,2}$ has rank $p^{2r}$.
Then $VD=\langle F^{r-1} \rangle \oplus \langle V, \ldots, V^{r-1} \rangle=N_r$.  
Also $V^2D= 0 \oplus \langle V^2, \ldots, V^{r-1} \rangle=N_{r-2}$. 
Thus $D$ has $a$-number $2$.
Continuing, one sees that $\nu_i=i-1$ for $1 \leq i \leq r-2$. 
In particular, $V$ is nilpotent on $D$ and thus the $p$-rank of $D$ is $0$.

More information on the final filtration is necessary to determine the final type of $I_{r,2}$.
First, $F^{-1}(N_r)=\langle F^{r-2}, F^{r-1} \rangle \oplus \langle 1, V, \ldots, V^{r-1} \rangle=N_{r+2}$.
Second, 
$F^{-r+2}(N_r)= \langle F, \ldots, F^{r-1} \rangle \oplus \langle 1, V, \ldots, V^{r-1} \rangle=N_{2r-1}$.
Then $VN_{2r-1}=0 \oplus \langle V, \ldots, V^{r-1} \rangle=N_{r-1}$ and
$VN_{r-1}= \langle V^2, \ldots, V^{r-1} \rangle=N_{r-2}$.
Thus $\nu_{r-1}=r-2$.
Then $I_{r,2}$ has final type $[0,1, \ldots, r-3,r-2,r-2]$ and Young type $\{r,1\}$. 
\end{proof}

\begin{lemma}
When $g=3$ or $g=4$, there is a unique indecomposable symmetric $BT_1$ group scheme 
$I_{g,2}$ of rank $p^{2g}$ with $p$-rank $0$ and $a$-number $2$.  
\end{lemma}

\begin{proof}
Let ${\mathbb G}$ be a symmetric $BT_1$ group scheme with rank $p^{2g}$, $p$-rank $0$ and $a$-number $2$.
Its Young type is $\{g,i\}$ for some $i \in \{1, \ldots, g-1\}$.
There are exactly $\lfloor g/2 \rfloor$ such group schemes which are decomposable, namely 
$I_{r,1} \oplus I_{g-r,1}$ for $1 \leq r \leq g/2$.
Thus there is a unique such ${\mathbb G}$ which is indecomposable when $g=3$ of $g=4$.
By Lemma \ref{a=2}, it is $I_{g,2}$.
\end{proof}

\subsection{One more indecomposable group scheme of dimension four}

There is one more indecomposable group scheme which occurs for dimension $g \leq 4$, 
which we denote by $I_{4,3}$.  It has covariant Dieudonn\'e module   
$D(I_{4,3})={\mathbb E}/(F^2-V)_\ell \oplus {\mathbb E}/(F-V)_\ell \oplus {\mathbb E}/(V^2-F)_\ell$.
Then $D(I_{4,3})$ has basis
$\langle 1,F,F^2 \rangle \oplus \langle 1,F \rangle \oplus \langle 1,V,V^2 \rangle$. 
One can check that $V^2D$ has basis $0 \oplus 0 \oplus \langle V^2 \rangle$ 
and thus $I_{4,3}$ has $a$-number $3$.
Also $I_{4,3}$ has $p$-rank $0$ since $V$ acts nilpotently on $D(I_{4,3})$.  
By the process of elimination, $I_{4,3}$ has final type $[0,0,1,1]$ and Young type $\{4,3,1\}$.

\section{Complete tables for dimension up to four} \label{Stables}

For convenience, we provide tables for dimension $g \leq 4$.
Some parts of these tables can be found in \cite{EV}.

\subsection{The case $g=1$:}

\renewcommand{\arraystretch}{1.065}
\[\begin{array}{|l|c|c|c|c|l|l|r|}
\hline
\mbox{Name} & \codim & f & a & \nu & \mu & \omega & \mbox{cycle class (reduced)} \\ 
\hline
\hline
L & 0 & 1 & 0 &  [1] & \emptyset & s_1 & \lambda_0 \\ 
\hline
I_{1,1} & 1 & 0 & 1 &  [0] & \{1\} & 1 & (p-1)\lambda_1 \\
\hline
\end{array}
\]

\subsection{The case $g=2$:}

\renewcommand{\arraystretch}{1.065}
\[\begin{array}{|l|c|c|c|c|l|l|r|}
\hline
\mbox{Name} & \codim & f & a &  \nu & \mu & \omega & \mbox{cycle class (reduced)} \\ 
\hline
\hline
L^2 & 0 & 2 & 0 &  [1,2] & \emptyset & s_2s_1s_2 & \lambda_0 \\
\hline
L \oplus I_{1,1} & 1 & 1 & 1 &  [1,1] & \{1\} & s_1s_2 & (p-1)\lambda_1 \\
\hline
I_{2,1} & 2 & 0 & 1 &  [0,1] & \{2\} & s_2 & (p-1)(p^2-1)\lambda_2 \\ 
\hline
I_{1,1}^2 & 3 & 0 & 2 &  [0,0] & \{2,1\} & 1 & (p-1)(p^2+1)\lambda_1\lambda_2 \\
\hline
\end{array}
\]

This is the smallest dimension for which the Newton polygon of $A$ does not determine the group scheme $A[p]$.
The Newton polygon $2G_{1,1}$ (supersingular, with four slopes of $1/2$)
occurs for both $(I_{1,1})^2$ and $I_{2,1}$.

\subsection{The case $g=3$:}


\renewcommand{\arraystretch}{1.065}
\[\begin{array}{|l|c|c|c|c|l|l|r|}
\hline
\mbox{Name} & \codim & f & a & \nu & \mu & \omega & \mbox{cycle class (reduced)} \\
\hline
\hline
L^3 & 0 & 3 & 0 &          [1,2,3] & \emptyset & s_3s_2s_3s_1s_2s_3 & \lambda_0\\
\hline
L^2 \oplus I_{1,1} & 1 & 2 & 1 & [1,2,2] & \{1\} & s_2s_3s_1s_2s_3 & (p-1)\lambda_1\\
\hline
L \oplus I_{2,1} & 2 & 1 & 1&    [1,1,2] & \{2\} & s_3s_1s_2s_3 & (p-1)(p^2-1)\lambda_2\\
\hline
L \oplus I_{1,1}^2 & 3 & 1 & 2 & [1,1,1] & \{2,1\} & s_1s_2s_3 & -(p-1)(p^2+1)\lambda_1\lambda_2 -2(p^3-1)\lambda_3\\
\hline
I_{3,1} & 3 & 0 & 1 &          [0,1,2] & \{3\} & s_3s_2s_3 & (p-1)(p^2-1)(p^3-1)\lambda_3\\
\hline
I_{3,2} & 4 & 0 & 2 &          [0,1,1] & \{3,1\} & s_2s_3 & (p-1)^2(p^2-p+1)\lambda_1\lambda_3\\
\hline
I_{1,1} \oplus I_{2,1} & 5 & 0 & 2 &   [0,0,1] & \{3,2\} & s_3 & -(p-1)^3(p+1)(p^2-p+1)(p^2+p+1)\lambda_2\lambda_3\\
\hline
I_{1,1}^3 & 6 & 0 & 3 &          [0,0,0] & \{3,2,1\} & 1 & (p-1)(p^2+1)(p^3-1)\lambda_1\lambda_2\lambda_3\\
\hline
\end{array}
\]

This is the smallest dimension for which the group scheme $A[p]$ does not determine the Newton polygon of $A$.
If $A[p] \simeq I_{3,1}$, then the Newton polygon of $A$ is usually $G_{1,2}+G_{2,1}$ 
(three slopes of $1/3$ and of $2/3$) 
but by \cite[Thm.\ 5.12]{O:hypsup} 
it can also be $3G_{1,1}$ (supersingular, with six slopes of $1/2$).

\subsection{The case $g=4$:}  


\renewcommand{\arraystretch}{1.065}
\[\begin{array}{|l|c|c|c|c|l|l|}
\hline
\mbox{Name} & \codim & f & a &  \nu & \mu & \omega\\ 
\hline
\hline
L^4                 & 0 & 4 & 0 & [1,2,3,4] & \emptyset    & 
s_4s_3s_4s_2s_3s_4s_1s_2s_3s_4 \\ \hline
L^3 \oplus I_{1,1}        & 1 & 3 & 1 & [1,2,3,3] & \{1\}     & 
s_3s_4s_2s_3s_4s_1s_2s_3s_4 \\ \hline
L^2 \oplus I_{2,1}        & 2 & 2 & 1 & [1,2,2,3] & \{2\}     & 
s_4s_2s_3s_4s_1s_2s_3s_4 \\ \hline
L^2 \oplus I_{1,1}^2      & 3 & 2 & 2 & [1,2,2,2] & \{2,1\}   & 
s_2s_3s_4s_1s_2s_3s_4 \\ \hline
L \oplus I_{3,1}        & 3 & 1 & 1 & [1,1,2,3] & \{3\}     & 
s_4s_3s_4s_1s_2s_3s_4 \\ \hline
L \oplus I_{3,2}        & 4 & 1 & 2 & [1,1,2,2] & \{3,1\}   & 
s_3s_4s_1s_2s_3s_4 \\ \hline
I_{4,1}                 & 4 & 0 & 1 & [0,1,2,3] & \{4\}     & 
s_4s_3s_4s_2s_3s_4 \\ \hline
L \oplus I_{1,1} \oplus I_{2,1} & 5 & 1 & 2 & [1,1,1,2] & \{3,2\}   & 
s_4s_1s_2s_3s_4\\ \hline
I_{4,2}                 & 5 & 0 & 2 & [0,1,2,2] & \{4,1\}   & 
s_3s_4s_2s_3s_4 \\ \hline
L \oplus I_{1,1}^3      & 6 & 1 & 3 & [1,1,1,1] & \{3,2,1\} & 
s_1s_2s_3s_4 \\ \hline
I_{1,1} \oplus I_{3,1}      & 6 & 0 & 2 & [0,1,1,2] & \{4,2\}   & 
s_4s_2s_3s_4 \\ \hline
I_{1,1} \oplus I_{3,2}      & 7 & 0 & 3 & [0,1,1,1] & \{4,2,1\} & 
s_2s_3s_4 \\ \hline
I_{2,1} \oplus I_{2,1}      & 7 & 0 & 2 & [0,0,1,2] & \{4,3\}   & 
s_4s_3s_4 \\ \hline
I_{4,3}                 & 8 & 0 & 3 & [0,0,1,1] & \{4,3,1\} & 
s_3s_4 \\ \hline
I_{1,1}^2 \oplus I_{2,1}    & 9 & 0 & 3 & [0,0,0,1] & \{4,3,2\} & 
s_4 \\ \hline
I_{1,1}^4               & 10& 0 & 4 & [0,0,0,0] &\{4,3,2,1\}& 
1 \\ \hline
\hline
\end{array}
\]

The cycle classes for this table can be found in \cite[15.3]{EV}.

It is not straight-forward to determine which Ekedahl-Oort strata lie in the boundary of which others.
When $g=4$, the answer to this question is given by the natural partial ordering on the Young type, 
which matches the Bruhat-Chevalley order on the elements of the Weyl group.
  
$$
\xygraph{
!{<0cm,0cm>;<1cm,0cm>;<0cm,1cm>::}
!{(0,5)*+{{\{\emptyset\}}}}="no"
!{(1,5)*+{{\{1\}}}}="1"
!{(2,5)*+{{\{2 \}}}}="2"
!{(3.2,5)*+{{\{3 \}}}}="3"
!{(4.6,5)*+{{\{4 \}}}}="4"
!{(2,4)*+{{\{2,1 \}}}}="21"
!{(3.2,4)*+{{\{3,1 \}}}}="31"
!{(4.6,4)*+{{\{4,1 \}}}}="41"
!{(3.2,3)*+{{\{3,2 \}}}}="32"
!{(4.6,3)*+{{\{4,2 \}}}}="42"
!{(6,3)*+{{\{4,3 \}}}}="43"
!{(3.2,2)*+{{\{3,2,1 \}}}}="321"
!{(4.6,2)*+{{\{4,2,1 \}}}}="421"
!{(6,2)*+{{\{4,3,1 \}}}}="431"
!{(7.4,2)*+{{\{4,3,2 \}}}}="432"
!{(9,2)*+{{\{4,3,2,1 \}}}}="4321"
"no":"1"
"1":"2"
"2":"3"
"3":"4"
"2":"21"
"3":"31"
"4":"41"
"21":"31"
"31":"41"
"31":"32"
"41":"42"
"32":"42"
"42":"43"
"32":"321"
"42":"421"
"43":"431"
"321":"421"
"421":"431"
"431":"432"
"432":"4321"
}$$

\bibliographystyle{alpha}
\bibliography{groupscheme}
\end{document}